\documentclass[reqno,12pt]{amsart}

\usepackage{mathrsfs}
\usepackage{amsmath,amssymb,slashed}
\usepackage{color,xcolor}
\usepackage[english]{babel}
\usepackage[utf8]{inputenc}
\usepackage{amsmath,amscd}
\usepackage{amsfonts}
\usepackage{amsthm}
\usepackage{graphicx}
\usepackage[colorinlistoftodos]{todonotes}
\usepackage{amsmath,amssymb,amsfonts}
\usepackage{underscore}
\usepackage[T1]{fontenc}
\usepackage[utf8]{inputenc}
\usepackage{amsthm}
\usepackage{mathrsfs}
\usepackage{amsmath,amssymb,slashed}
\usepackage{color,xcolor}

\input xy
\xyoption{all}

\theoremstyle{plain}
\newtheorem{thm}{Theorem}[section]
\newtheorem{prop}[thm]{Proposition}
\newtheorem{lem}[thm]{Lemma}
\newtheorem{cor}[thm]{Corollary}

\newtheorem{pro}[thm]{Problem}

\theoremstyle{definition}
\newtheorem{defn}[thm]{Definition}

\theoremstyle{remark}
\newtheorem{rmk}[thm]{Remark}
\newtheorem{ques}[thm]{Question}

\newcommand{\la}{\langle}
\newcommand{\ra}{\rangle}

\begin{document}

\title[The Takai duality of $L^p$-operator algebras]{The Takai duality of $L^p$-operator algebras with incompressibility}

\author[Z. Wang]{Zhen Wang}
\address{School of Mathematics\\Hangzhou Normal University\\Hangzhou 311121\\ P.~R. China}
\email{wangzhen@hznu.edu.cn}

\author[S. Zhu]{Sen Zhu}
\address{Department of Mathematics\\Jilin University\\Changchun 130012\\P.~R. China}
\email{zhusen@jlu.edu.cn}
\thanks{This work was partially supported by the National Natural Science Foundation of China (Grant numbers 12201240, 12171195) and a grant from Hangzhou
Normal University (Grant number 4235C50224204070).}


%\subjclass[2010]{Primary 22D35 , 47L65; Secondary 43A25, 47L10}

\begin{abstract}
Let $G$ be a countable discrete Abelian group, $A$ be a separable unital $L^p$-operator algebra which has unique $L^p$-operator matrix norms for $p\in [1,\infty)$, and $\alpha$ be an isometric action of $G$ on $A$. We prove in this paper that if $A$ is $p$-isometrically incompressible, then $F^{p}(\hat{G},F^p(G,A,\alpha),\hat{\alpha})$, the iterated $L^p$-operator crossed product, is isometrically isomorphic to $\overline{M}_{G}^{p}\otimes_{p}A$ if and only if $p=2$. Furthermore, we prove that if $A=M_n^p$, then $F^{p}(\hat{G},F^p(G,A,\alpha),\hat{\alpha})$ is isomorphic to $\overline{M}_{G}^{p}\otimes_{p}A$ if and only if either $p=2$ or $G$ is finite. This shows that the Takai duality for $C^*$-algebras
can not be generalized to $L^p$-operator algebras with $p\neq 2$, and solves a problem raised by N. C. Phillips in the negative.
\end{abstract}

\date{\today}
\maketitle

 %\tableofcontents

\section{Introduction}

This paper is a continuation of a recent paper \cite{WZ}, aiming to study the $L^p$-Takai duality problem raised by N. C. Phillips \cite[Problem 8.7]{PlpsOpenQues}. Before proceeding, we first recall some terminology and basic facts.

%The aim of this paper is to study the Takai duality problem for $L^p$ operator crossed products posed by N. C. Phillips \cite[Problem 8.7]{PlpsOpenQues}. To proceed, we
%first introduce some terminology.

A Banach algebra $A$ is called an $L^{p}$-{\it operator algebra} if it can be isometrically represented on an $L^{p}$-space $(p\in [1, \infty))$.
The study of $L^{p}$-operator algebras dates back to Herz's influential paper \cite{Herz} on harmonic analysis on $L^p$-spaces.
People's renewed interests in $L^p$-operator algebras were inspired by N. C. Phillips
\cite{Phi Cuntz, N. C. Phillips Lp}, who showed that several important constructions of $C^*$-algebras can be generalized to the setting of $L^p$-operator algebras.
Recently the study of $L^p$-operator algebras has received much attention, and
a number of authors have made significant contributions to it (see \cite{Blecher and Phillips,Braga,Choi and Gardela,Chung lp Roe 1,Cort,Gardella and Thiel groupoid,Gardella and Thiel,Gardella and Thiel quotient,Gardella and Thiel convolution operators,Gardella and Thiel convolution,Gardella modern,Hejazian and Pooya,N. C. Phillips Lp,WW,WZ,nuclear,W}).

The Takai duality theorem is a fundamental result concerning crossed products of $C^*$-algebras. Let $(G,A,\alpha)$ be a $C^{*}$-dynamical system, where $A$ is a $C^{*}$-algebra, $G$ is a locally compact group and $\alpha$ is a continuous homomorphism from $G$ to ${\rm Aut}(A)$ (the group of automorphisms of $A$).
If $G$ is Abelian, then there exists a dual system $(\hat{G},A\rtimes_{\alpha} G, \hat{\alpha})$.
The Takai duality theorem \cite{Raeburn,Takai} says that the iterated crossed product $(A\rtimes_{\alpha}G)\rtimes _{\hat{\alpha}}\hat{G}$ is $*$-isomorphic to $\mathcal{K}\left(L^{2}(G)\right)\otimes A$, where $\mathcal{K}\left(L^{2}(G)\right)$ denotes the algebra of compact operators on $L^{2}(G)$.  We refer to \cite{Katsoulis Takai,Katsoulis} for the Takai duality theorem for crossed products of non-selfadjoint operator algebras.

N. C. Phillips \cite{N. C. Phillips Lp} introduced the $L^p$-operator crossed products with the aim to compute the $K$-theory groups of $L^{p}$-Cuntz algebras.
Furthermore, he posed the following Takai duality problem for $L^p$-operator crossed products.

\begin{pro}[{\cite[Problem 8.7]{PlpsOpenQues}}]\label{P:main}
Let $p\in[1,\infty)$. Let $\alpha:G\rightarrow \mathrm{Aut}(A)$ be an isometric action of locally compact Abelian group on an $L^p$-operator algebra. Then there is dual action $$\hat{\alpha}:\hat{G}\rightarrow \mathrm{Aut}(F^p(G,A,\alpha)),$$
where $F^p(G,A,\alpha)$ is the full $L^p$-operator crossed product.
Is there an analog of Takai duality for the crossed products by this action?
\end{pro}

Inspired by the preceding problem, we are interested in the following question.

%it is natural to ask whether $F^{p}(\hat{G},F^p(G,A,\alpha),\hat{\alpha})$  and $\overline{M}_{G}^{p}\otimes_{p}A$ are isometrically isomorphic or simply isomorphic (see Questions \ref{Q:main-1} and \ref{Q:main-2}).

%In view of the Takai duality theorem for $C^*$-crossed products, the preceding question lead us to consider the following question.

\begin{ques}\label{Q:main-1}
Let $G, A$, $\alpha$ and $p$ be as in Problem \ref{P:main}. Is the iterated $L^p$-operator crossed product $F^{p}(\hat{G},F^p(G,A,\alpha),\hat{\alpha})$ isometrically isomorphic to $\overline{M}_{G}^{p}\otimes_{p}A$?
\end{ques}

We recall the definition of $\overline{M}_{G}^{p}$. Let $G$ be a countable discrete group.
Given a finite subset $F$ of $G$, we define $M_F^p$ to be the set of all $a\in \mathcal{B}(l^p(G))$ such that $a\xi=0$ whenever $\xi|_F=0$ and such that $a\xi\in l^p(F)\subset l^p(G)$ for $\xi\in l^p(F)$.
Then we define $$M_G^p=\bigcup_{F\subset G ~\mathrm{finite}}M_F^p $$
and define $\overline{M}_{G}^{p}$ to be the closure of $M_G^p$ in the operator norm on $\mathcal{B}(l^p(G))$ (see \cite[Example 1.6]{N. C. Phillips Lp}).
Note that $\overline{M}_{G}^{p}=\mathcal{K}(l^p(G))$ when $p>1$, and $\overline{M}_{G}^{p}$ is strictly smaller than $\mathcal{K}(l^1(G))$ when $p=1$ and $G$ is infinite (see \cite[Corollary 1.9 and Example 1.10]{N. C. Phillips Lp}).

%{\color{red}Here $\overline{M}_{G}^{p}$ is taken from Example 1.6 in \cite{N. C. Phillips Lp}. The precise definition will be provided in the end of subsection \ref{compact}.}

The first result of this paper is the following theorem, which gives a complete answer to Question \ref{Q:main-1} for $L^p$-operator algebras
with $p$-isometrical incompressibility.

\begin{thm}\label{MainThm}
Let $p\in [1,\infty)$, $G$ be a countable discrete Abelian group, and $A$ be a separable unital $L^p$-operator algebra with unique $L^p$-operator matrix norms. If $A$ is $p$-isometrically incompressible, then $F^{p}(\hat{G},F^{p}(G,A,\alpha),\hat{\alpha})$ is isometrically isomorphic to $\overline{M}_{G}^{p}\otimes_{p}A$ if and only if $p=2$.
\end{thm}

Recall that an $L^p$-operator algebra $A$ is said to be {\it $p$-isometrically incompressible} if for every $L^p$-space $E$, every contractive, injective homomorphism $\rho:A\rightarrow \mathcal{B}(E)$ is isometric (see {\cite[Definition 4.1]{Gardella and Thiel}}).

The $L^p$-matrix algebra $M_n^p$, $L^p$-Cuntz algebra $\mathcal{O}_d^p$,  and $C^*$-algebras are examples of $p$-isometrically incompressible $L^p$-operator algebras (see \cite[Theorem 7.2,~Corollary 8.10]{Phi Cuntz} and \cite[Theorem 10]{Bonsall}). The direct limit of $p$-isometrically incompressible $L^p$-operator algebras is still $p$-isometrically incompressible (see \cite[Lemma 4.2]{Gardella and Thiel}).

Our proof of Theorem \ref{MainThm} depends on a characterization of the $p$-isometric incompressibility of $F^p(\hat{G},F^p(G,A,\alpha),\hat{\alpha})$ and $\overline{M}_{G}^{p}\otimes_{p}A$ (see Proposition \ref{P:key}). The reader will see that the $p$-isometric incompressibility provides a key tool to determine whether $F^p(\hat{G},F^p(G,A,\alpha),\hat{\alpha})$ and $\overline{M}_{G}^{p}\otimes_{p}A$
are isometrically isomorphic.

%shows that the
%$p$-isometrical incompressibility provides important the key inv
%
% (see Proposition \ref{P:key}).

Note that the isometric isomorphism and the isomorphism for Banach algebras do not coincide.
Theorem \ref{MainThm} and Problem \ref{P:main} inspire the following natural question.

\begin{ques}\label{Q:main-2}
Let $G, A$, $\alpha$ and $p$ be as in Problem \ref{P:main}. Is the iterated $L^p$-operator crossed product $F^{p}(\hat{G},F^p(G,A,\alpha),\hat{\alpha})$ isomorphic to $\overline{M}_{G}^{p}\otimes_{p}A$?
\end{ques}

Inspired by D. Williams' proof (see \cite[Theorem 7.1]{Williams}) for the Takai duality theorem for crossed products of $C^*$-algebras, we constructed in \cite{WZ} a homomorphism $\Phi$ from the iterated $L^p$-operator crossed product $F^{p}(\hat{G},F^p(G,A,\alpha),\hat{\alpha})$ to $\overline{M}_{G}^{p}\otimes_{p}A$, which is a natural $L^p$-analog of D. Williams' map. For countable discrete Abelian groups $G$ and separable unital $L^p$-operator algebras $A$ that have unique $L^p$-operator matrix norms, we prove that $\Phi$ is an isomorphism if and only if either $G$ is finite or $p=2$ (see \cite[Theorem 1.2]{WZ}).  This suggests a negative answer to Question \ref{Q:main-2} for the case $p\ne 2$.

The second result of this paper is the following theorem, which answers Question \ref{Q:main-1} for the $L^p$-operator algebra $A=M_n^p$.

\begin{thm}\label{MainThm2}
Let $p\in [1,\infty)$, $G$ be a countable discrete Abelian group, and $A=M_n^p$. Then $F^{p}(\hat{G},F^{p}(G,A,\alpha),\hat{\alpha})$ is isomorphic to $\overline{M}_{G}^{p}\otimes_{p}A$ if and only if $p=2$ or $G$ is finite.
\end{thm}

The proof of Theorem \ref{MainThm2} follows a similar line to that of Theorem \ref{MainThm}, and depends on a characterization of the $p$-incompressibility of $F^p(\hat{G},F^p(G,A,\alpha),\hat{\alpha})$ (see Proposition \ref{P:key2}).
The reader will see that the $p$-incompressibility provides a key tool to determine whether $F^p(\hat{G},F^p(G,A,\alpha),\hat{\alpha})$ and $\overline{M}_{G}^{p}\otimes_{p}A$ isomorphic in the case that $A=M_n^p$.

\begin{rmk}
From Theorems \ref{MainThm} and \ref{MainThm2}, one can see that the Takai duality for $C^*$-algebras
can not be generalized to $L^p$-operator algebras with $p\ne 2$.
This gives a negative answer to Problem \ref{P:main}. The fundamental reason why Takai duality cannot be extended to general $L^p$-operator algebras with $p\neq 2$ is that, for infinite groups $G$, the Fourier transform $$\mathcal{F}_p:L^p(G)\rightarrow L^{q}(\hat{G})$$ is not an isomorphism (see \cite[page 227]{Hew}), where $q$ is the conjugate exponent of $p$, that is, $1/p+1/q=1$. When $p=2$, the Fourier transform $\mathcal{F}_2$ is an isometric isomorphism. This result is known as Plancherel's theorem and is very useful for proving Takai duality for $C^*$-algebras (see \cite[7.9.1]{Pedersen}).
\end{rmk}

% It implies that the Gelfand transform $$\Gamma_2:C^*_\lambda(\hat{G})\rightarrow C_0(G)$$ is an isometric isomorphism, since $$\Gamma_2(f)=\mathcal{F}_2^{-1}f\mathcal{F}_2$$ for all $f\in C_c(\hat{G})$, where $C^*_\lambda(\hat{G})$ is the reduced group $C^*$-algebra of $\hat{G}$.
%
%
%
%However, this cannot hold for the reduced group $L^p$-operator algebra $F^p_\lambda(\hat{G})$ when $p\neq 2$. Actually, the Gelfand transform of $F^p_\lambda(\hat{G})$ has similar properties to $L^p$-Fourier transform, that is, the Gelfand transform $$\Gamma_p:F^p_\lambda(\hat{G})\rightarrow C_0(G)$$ is an isomorphism if and only if $p=2$ or $G$ is finite (see \cite[Corollary 3.20]{Gardella and Thiel}).
%
%

The following is an immediate corollary of  Theorems \ref{MainThm} and \ref{MainThm2}.

\begin{cor}
Let $p\in [1,\infty)$, $G$ be a countably infinite discrete Abelian group, and $A=M_n^p$.
Then the following are equivalent:
\begin{enumerate}
\item[(i)] $F^{p}(\hat{G},F^{p}(G,A,\alpha),\hat{\alpha})$ is isometrically isomorphic to $\overline{M}_{G}^{p}\otimes_{p}A$;
\item[(ii)] $F^{p}(\hat{G},F^{p}(G,A,\alpha),\hat{\alpha})$ is isomorphic to $\overline{M}_{G}^{p}\otimes_{p}A$;
\item[(iii)] $p=2$.
\end{enumerate}
\end{cor}

The rest of this paper is organized as follows.
In Section 2, we shall make some preparation.
Section 3 is devoted to the proofs of Theorems \ref{MainThm} and \ref{MainThm2}.

\section{Preliminaries}

%In this section, we recall some terminology and notations.

In this section, we review the terminology and notation used in the sequel.

%make some preparations for the proof of Theorem \ref{MainThm}.

\subsection{$L^p$-operator crossed products} \label{compact}
Recall that a Banach algebra $A$ is an $L^{p}$-{\it operator algebra} if it can be isometrically represented on an $L^{p}$-space $(p\in [1, \infty))$.
Let $\pi:A\rightarrow \mathcal{B}(L^p(X,\nu))$ be a representation of $A$, where $\mathcal{B}(L^p(X,\nu))$
denotes the algebra of all bounded linear operators on $L^p(X,\nu)$. We say that $\pi$ is {\it $\sigma$-finite} if $\nu$ is $\sigma$-finite.

Let $G$ be a locally compact group. Then there exists a left Haar measure $\mu$ on $G$. An $L^{p}$-{\it operator algebra dynamical system} is a triple $(G,A,\alpha)$ consisting of a locally compact group $G$, an $L^{p}$-operator algebra $A$ and a continuous homomorphism $\alpha:G\rightarrow \mathrm{Aut}(A)$, where $\mathrm{Aut}(A)$ is the group of isometric automorphisms of $A$.
A {\it contractive covariant representation} of $(G,A,\alpha)$ on an $L^{p}$-space $E$ is a pair $(\pi,v)$ consisting of a nondegenerate contractive homomorphism $\pi:A\rightarrow \mathcal{B}(E)$ and an isometric group representation $v:G\rightarrow \mathcal{B}(E)$, satisfying the covariance condition $$v_{t}\pi(a)v_{t^{-1}}=\pi(\alpha_{t}(a))$$ for $t\in G$ and $a\in A$.
A covariant representation $(\pi,v)$ of $(G,A,\alpha)$ is {\it $\sigma$-finite} if $\pi$ is $\sigma$-finite.

Denote by $C_{c}(G,A,\alpha)$ the vector space of continuous compactly supported functions $G\rightarrow A$, made into an algebra over $\mathbb{C}$ with product given by twisted convolution, that is, $$(f*g)(t):=\int_{G}f(s)\alpha_{s}(g(s^{-1}t))d\mu(s)$$ for $f,g\in C_{c}(G,A,\alpha)$ and $t\in G$. The {\it integrated form} of $(\pi,v)$ is the nondegenerate contractive homomorphism $\pi\rtimes v: C_{c}(G,A,\alpha)\rightarrow \mathcal{B}(E)$ given by $$(\pi\rtimes v)(f)(\xi):=\int_{G}\pi(f(t))v_{t}(\xi)d\mu(t)$$ for $f\in C_{c}(G,A,\alpha)$ and $\xi\in E$. Denote by $\mathrm{Rep}_{p}(G,A,\alpha)$ the class of all nondegenerate $\sigma$-finite contractive covariant representations of $(G,A,\alpha)$ on $L^{p}$-spaces. The {\it full $L^p$-operator crossed product} $F^{p}(G,A,\alpha)$ is defined as the completion of $C_{c}(G,A,\alpha)$ in the norm $$||f||_{F^{p}(G,A,\alpha)}:=\sup\big\{||(\pi\rtimes v)(f)||:(\pi,v)\in \mathrm{Rep}_{p}(G,A,\alpha)\big\}.$$

Let $(G,A,\alpha)$ be an $L^{p}$-operator algebra dynamical system.
Given a nondegenerate $\sigma$-finite contractive representation $\pi_{0}:A\rightarrow \mathcal{B}(E_{0})$ on an $L^{p}$-space $E_{0}$, its associated {\it regular covariant representation} is the pair $(\pi,\lambda_{p}^{E_{0}})$ on $L^{p}(G)\otimes_{p} E_{0}\cong L^{p}(G,E_{0})$ given by
$$\pi(a)(\xi)(s):=\pi_{0}\left(\alpha_{s^{-1}}(a)\right)(\xi(s)) \quad$$ and $$\lambda_p^{E_0}(s)(\xi)(t):=\xi(s^{-1}t)$$
for $a\in A$, $\xi\in L^{p}(G,E_{0})$, and $s,t\in G$.
We denote by $\mathrm{RegRep}_{p}(G,A,\alpha)$ the class consisting
of nondegenerate $\sigma$-finite contractive regular covariant representations of $(G,A,\alpha)$, which is a subclass of $\mathrm{Rep}_{p}(G,A,\alpha)$. The {\it reduced $L^p$-operator crossed product} $F^{p}_{\lambda}(G,A,\alpha)$ is defined as the completion of $C_{c}(G,A,\alpha)$ in the norm $$||f||_{F^{p}_{\lambda}(G,A,\alpha)}:=\sup\{||(\pi\rtimes v)(f)||:(\pi,v)\in \mathrm{RegRep}_{p}(G,A,\alpha)\}.$$

By \cite[Theorem 7.1]{Phillips look like}, if $G$ is amenable, then the identity map on $C_c(G,A,\alpha)$ can be extended to an isometric isomorphism from $F^{p}(G,A,\alpha)$ onto $ F^{p}_{\lambda}(G,A,\alpha)$.
If $A=\mathbb{C}$, then it is easy to see that $F^{p}(G,A,\mathrm{id})$ is the full group $L^{p}$-operator algebra $F^{p}(G)$ and $F^{p}_{\lambda}(G,A,\mathrm{id})$ is the reduced group $L^{p}$-operator algebra $F^{p}_{\lambda}(G)$,
where $\mathrm{id}$ is the trivial action of $G$ on $\mathbb{C}$.

If, in addition, $G$ is Abelian, then we let $\hat{G}$ denote the dual group of $G$. For each $\gamma\in \hat{G}$, N. C. Phillips defined an isomorphism $\hat{\alpha}_{\gamma}:C_{c}(G,A,\alpha)\rightarrow C_{c}(G,A,\alpha)$ by $\hat{\alpha}_{\gamma}(f)(s):=\overline{\gamma(s)}f(s)$ for $f\in C_{c}(G,A,\alpha)$ and $s\in G$ (see \cite[Definition 3.15]{N. C. Phillips Lp}). Thus $\hat{\alpha}_{\gamma}$ can extend to an isometry on $F^{p}(G,A,\alpha)$ by continuity (see \cite[Theorem 3.18]{N. C. Phillips Lp}).
Hence there is a {\it dual system} $(\hat{G},F^{p}(G,A,\alpha),\hat{\alpha})$.
So we obtain the iterated $L^{p}$-operator crossed product $F^{p}(\hat{G},F^{p}(G,A,\alpha),\hat{\alpha})$.

\subsection{Unique $L^p$-operator matrix norms}
In this subsection, we recall the definition of unique $L^p$-operator matrix norms.
For each positive integer $n$, let $l_n^p=l^p\left(\{1,2,\cdots,n\},\nu\right)$, where $\nu$ is the counting measure on $\{1,2,\cdots,n\}$.
We denote $M_n^p=\mathcal{B}(l_n^p)$.
\begin{defn}\label{matrix algebra}
Given a closed subalgebra $A$ of $\mathcal{B}(L^p(X,\mu))$, we denote by $M_n^p\otimes_p A$ the {\it $L^p$-matrix algebra}, that is, the Banach subalgebra of $\mathcal{B}(L^p(\{1,2,\cdots,n\}\times X,\nu\times\mu))$ generated by all $T\otimes a$ for $T\in M_n^p$ and $a\in \mathcal{B}(L^p(X,\mu))$. Clearly, each element of $M_n^p\otimes_p A$ is of form $[a_{i,j}]_{1\leq i,j\leq n}$ with $a_{i,j}\in A$,
which is also written as $\sum_{i,j=1}^n e_{i,j}\otimes a_{i,j}$, where $\{e_{i,j}\}_{1\leq i,j\leq n}$ are the canonical matrix units of $M_n^p$.
\end{defn}

\begin{defn} %\cite{Lp AF}
Let $A$ be a closed subalgebra of $\mathcal{B}(L^p(X,\mu))$, $B$ be a closed subalgebra of $\mathcal{B}(L^p(Y,\nu))$
and $\varphi$ be a linear map $\varphi: A\rightarrow B$.  We denote by $\mathrm{id}_{M_n^p}\otimes \varphi$ the map from $M_n^p\otimes_p A$ to $M_n^p\otimes_p B$ defined by $$\mathrm{id}_{M_n^p}\otimes \varphi\left(\sum_{i,j=1}^n e_{i,j}\otimes a_{i,j}\right)=\sum_{i,j=1}^n e_{i,j}\otimes \varphi(a_{i,j}) $$ for
$\sum_{i,j=1}^n e_{i,j}\otimes a_{i,j}\in M_n^p\otimes_p A$.
We denote $\|\varphi\|_{cb}=\sup_{n\in \mathbb{Z}_{>0}}\|\mathrm{id}_{M_n^p}\otimes \varphi\|$.
We say that $\varphi$ is {\it $p$-completely contractive} if $\|\varphi\|_{cb}\leq 1$, and say that $\varphi$ is {\it $p$-completely isometric} if  $\mathrm{id}_{M_n^p}\otimes \varphi$ is isometric for all positive integer $n$.
\end{defn}

\begin{defn}[{\cite[Definition 4.1]{Lp AF}}]\label{matrix norm}
Let $p\in [1,\infty)$ and $A$ be a separable $L^p$-operator algebra. We say that $A$ has
{\it unique $L^p$-operator matrix norms} if whenever $(X,\mathcal{B},\mu)$ and $(Y,\mathcal{C},\nu)$ are $\sigma$-finite measure spaces such that $L^p(X,\mu)$ and $L^p(Y,\nu)$ are separable, $\pi:A\rightarrow \mathcal{B}(L^p(X,\mu))$ and $\sigma:A\rightarrow \mathcal{B}(L^p(Y,\nu))$ are isometric representations, and $\pi(A)$ and $\sigma(A)$ are given the matrix normed structure of Definition \ref{matrix algebra}, then $\sigma\circ\pi^{-1}:\pi(A)\rightarrow \sigma(A)$ is completely isometric.
\end{defn}

$M_n^p$ and $C(X)$ are basic examples with unique $L^p$-operator matrix norms (see \cite[Corollary 4.4 \& Proposition 4.6]{Lp AF}), where $X$ is a compact metrizable space.

\subsection{$C^*$-cores}
The $C^*$-core of an $L^p$-operator algebra was introduced by Y. Choi,  E. Gardella and H. Thiel \cite{Choi and Gardela}, which is an invariant under isometric isomorphisms of $L^p$-operator algebras.

\begin{defn}[{\cite[Definition 2.6.1]{Palmer}}]
Let $A$ be a unital Banach algebra in which $\|1\|=1$. The {\it numerical range} $W(a)$ of an element $a$ in $A$ is the set of all numbers $w(a)\in \mathbb{C}$ with $w$ being a continuous linear functional on $A$ with $\|w\|=w(1)=1$.
\end{defn}

\begin{defn}[{\cite[Definition 5.5]{Lp AF}}]
Let $A$ be a unital Banach algebra in which $\|1\|=1$.
An element $a\in A$ is said to be {\it hermitian} if $W(a)\subset \mathbb{R}$.
\end{defn}

\begin{defn}[{\cite[Definition 2.10]{Choi and Gardela}}]
Let $p\in [1,\infty)$, and let $A$ be a unital $L^p$-operator algebra. Denote by $A_h$ the set of hermitian elements in $A$. The algebra $\mathrm{core}(A):={A_h}+\mathrm{i}{A_h}$ is called the $C^*$-{\it core} of $A$.
\end{defn}

By \cite[Theorem 2.9]{Choi and Gardela}, the $C^*$-core of $A$ is the largest unital $C^*$-subalgebra of $A$, and it is commutative when $p\in [1,\infty)\setminus\{2\}$.
By \cite[Proposition 2.13]{Choi and Gardela}, if $\varphi:A\rightarrow B$ is an isometric isomorphism between two $L^p$-operator algebras $A$ and $B$, then $\varphi: \mathrm{core}(A)\rightarrow \mathrm{core}(B)$ is a $*$-isomorphism. This implies that the $C^*$-core is invariant under isometric isomorphism between $L^p$-operator algebras.

\subsection{$p$-incompressibility}\label{in}
To determine whether the iterated $L^p$-operator crossed product $F^{p}(\hat{G},F^p(G,A,\alpha),\hat{\alpha})$ and $\overline{M}_{G}^{p}\otimes_p A$ are isomorphic, we introduce $p$-incompressible $L^p$-operator algebras. This property is invariant under isomorphisms of $L^p$-operator algebras. The definition is inspired by the work of N. C. Phillips and B. Johnson on incompressible Banach algebras \cite{incompressible}.

\begin{defn}
Let $A$ be an $L^p$-operator algebra.
\begin{enumerate}
\item[(i)]  We say that $A$ is {\it $p$-incompressible} if whenever  $E$ is an $L^p$-space and $\varphi:A\rightarrow \mathcal{B}(E)$ is a bounded injective homomorphism, then $\varphi$ is bounded below.
\item[(ii)] We say that $A$ is {\it $p$-uniformly incompressible} if for every $R\in [1,\infty)$ there is $\epsilon>0$ such that whenever $E$ is an $L^p$-space and $\varphi:A\rightarrow \mathcal{B}(E)$ is an injective homomorphism with $\|\varphi\|\leq R$, then for each $a\in A$ we have $\|\varphi(a)\|\geq \epsilon \|a\|$.
\end{enumerate}
\end{defn}

\begin{rmk}
\begin{enumerate}
\item[(i)] One can see that $p$-uniform incompressibility implies $p$-incompressibility.
 $C^*$-algebras and $M_n^p$ are examples of $p$-uniformly incompressible $L^p$-operator algebras (see \cite{incompressible}).
 \item[(ii)] Let $A$ and $B$ be two $L^p$-operator algebras. If $A$ and $B$ are isomorphic, then $B$ is $p$-incompressible if and only if $A$ is $p$-incompressible.
\end{enumerate}
\end{rmk}
An example of $L^p$-operator algebra that is neither $p$-incompressible nor $p$-isometriccally incompressible is the reduced group $L^p$-operator algebra $F^p_\lambda(G)$ when $G$ is infinite Abelian and $p\in [1,\infty)\setminus\{2\}$ (see Theorem \ref{Gelfand} and Lemma \ref{below}).
%The reason follows from the Gelfand transform of $F^p_\lambda(G)$, which differs from the case $p=2$.

\begin{thm}[{\cite[Corollary 3.20]{Gardella and Thiel}}]\label{Gelfand}
Let $p\in [1,\infty)\setminus \{2\}$. Let $G$ be a locally compact Abelian group. Then  Gelfand transform $\Gamma_p: F^p_\lambda(G)\rightarrow C_0(\hat{G})$ is contractive injective with dense range and $\Gamma_p$ is surjective if and only if $G$ is finite or $p=2$.
\end{thm}

\begin{lem}[{\cite[Exercise 5 on page 93]{Conway}}]\label{below}
Let $X$ and $Y$ be Banach spaces and let $T\in \mathcal{B}(X,Y)$. Then there is a constant $c>0$ such that $\|Tx\||\geq c\|x\|$ for all $x\in X$ if and only if $\ker T=\{0\}$ and $\mathrm{ran}~T$ is closed.
\end{lem}

Theorem \ref{Gelfand} and Lemma \ref{below} show that $F^p_\lambda(G)$ is neither $p$-isometrically incompressible nor $p$-incompressible when $G$ is infinite and $p\neq 2$. This is the key reason why the Takai duality for $C^*$-algebras can not be generalized to the $L^p$-case in general.

Since $C_0(\hat{G})$ is both $p$-isometrically incompressible and $p$-incompressible, we have the following corollary.

\begin{cor}\label{1}
Let $p\in [1,\infty)$ and $G$ be an infinite locally compact Abelian group. Then the following are equivalent:
\begin{enumerate}
\item[(i)] $F^p_\lambda(G)$ is isometrically isomorphic to $C_0(\hat{G})$;
\item[(ii)] $F^p_\lambda(G)$ is isomorphic to $C_0(\hat{G})$;
\item[(iii)] $p=2$.
\end{enumerate}
Moreover, if $G$ is finite, then $F^p_\lambda(G)$ is always isomorphic to $C(\hat{G})$ and they are isometrically isomorphic if and only if $p=2$ or $G$ is trivial.
\end{cor}
\begin{proof}
The implications $(iii)\Rightarrow (i) \Rightarrow (ii)$ are obvious. Now we prove $(ii)\Rightarrow (iii)$.

If $p\neq 2$ and $G$ is infinite, then $F^p_\lambda(G)$ is not $p$-incompressible. However, $C_0(\hat{G})$ is $p$-incompressible and isomorphisms between $L^p$-operator algebras preserve $p$-incompressibility, a contradiction. Thus $p=2$.

Now if $G$ is a finite group, from the Gelfand transform of $F^p_\lambda(G)$, we have that $F^p_\lambda(G)$ is always isomorphic to $C(\hat{G})$.
By \cite[Corollary 2.21]{Choi and Gardela}, it follows that
$\mathrm{core}(F^p_\lambda(G))=\mathbb{C}$ and $\mathrm{core}(C(\hat{G}))=C(\hat{G})$.
Since isometric isomorphisms between $L^p$-operator algebras preserve $C^*$-cores, it follows that $F^p_\lambda(G)$ is isometrically isomorphic to $C(\hat{G})$ if and only if $p=2$ or $G$ is trivial.
\end{proof}

Since the Takai duality relies on the Gelfand transform, we will see that the Takai duality for $L^p$-operator crossed products share similar properties of Corollary \ref{1}.

\section{Proof of Theorem \ref{MainThm}}

To prove our main Theorem, we recall the main result of our previous paper \cite{WZ}.
In \cite{WZ}, we construct $L^p$-analogues of the construction (\cite[Theorem 7.1]{Williams}) as follows:
%$$\begin{CD}
%F^{p}(\hat{G},F^{p}(G,A,\alpha),\hat{\alpha}) @>\Phi_1>> F^{p}(G,
%F^{p}(\hat{G},A,\beta),\hat{\beta}\otimes \alpha) @ >\Phi_2>>F^{p}(G,C_{0}(G,A),\mathrm{lt}\otimes\alpha) \\
%@.     @. @VV\Phi_3V\\
%@.      \overline{M}_{G}^{p}\otimes_{p}A @<\Phi_4<< F^{p}(G,C_{0}(G,A),\mathrm{lt}\otimes\mathrm{id})
%\end{CD} .$$

$$\begin{CD}
F^{p}\left(\hat{G},F^{p}(G,A,\alpha),\hat{\alpha}\right) @>\Phi_1>> F^{p}\left(G,
F^{p}(\hat{G},A,\beta),\hat{\beta}\otimes \alpha\right)\\
@.   @VV\Phi_2V\\
F^{p}\left(G,C_{0}(G,A),\mathrm{lt}\otimes\mathrm{id}\right) @<\Phi_3<<F^{p}\left(G,C_{0}(G,A),\mathrm{lt}\otimes\alpha\right)\\
@VV\Phi_4V @. \\
\overline{M}_{G}^{p}\otimes_{p}A
\end{CD}$$

The main result of \cite{WZ} is the following theorem.

\begin{thm}[{\cite[Theorem 1.2]{WZ}}]\label{WZ}
Let $(G,A,\alpha)$ be an $L^p$-operator algebra dynamical system,
where $G$ is a countable discrete Abelian group, and $A$ is a separable unital $L^p$-operator algebra which has unique $L^p$-operator matrix norms.
Let $\Phi=\Phi_4\circ\Phi_3\circ\Phi_2\circ\Phi_1$.
 Then
\begin{enumerate}
\item[(i)] $\Phi_{1}$, $\Phi_{3}$ and $\Phi_4$ are three isometric isomorphisms for $p\in [1,\infty)$;
\item[(ii)] $\Phi_2$ is an isomorphism if and only if either $G$ is finite or $p=2$; in particular, $F^2(\hat{G},F^2(G,A,\alpha),\hat{\alpha})$ is isometrically isomorphic to $F^{p}(G,C_{0}(G,A),\mathrm{lt}\otimes\alpha)$;
\item[(iii)] $\Phi:F^p(\hat{G},F^p(G,A,\alpha),\hat{\alpha})\rightarrow  \overline{M}_{G}^{p}\otimes_p A$ is an isomorphism if and only if either $G$ is finite or $p=2$; in particular, $F^2(\hat{G},F^2(G,A,\alpha),\hat{\alpha})$ is isometrically isomorphic to $ \overline{M}_{G}^{2}\otimes_2 A$.
\end{enumerate}
\end{thm}

\begin{rmk}\label{2}
Observe that $\Phi_2$ is induced by the Gelfand transform $\Gamma_p:F^p_\lambda(\hat{G})\rightarrow C_0(G)$. The map $\Gamma_p$ is injective and contractive, has dense range, and is not surjective when $p\neq 2$ and $G$ is infinite. Consequently, $\Phi_2$ inherits all these properties. This fact plays an essential role in the proofs of Propositions \ref{P:key} and \ref{P:key2}.
\end{rmk}

We now prove Theorem \ref{MainThm}. This theorem is a direct consequence of the following proposition, since isometric isomorphisms between
$L^p$-operator algebras preserve the $p$-isometric incompressibility.

\begin{prop}\label{P:key}
Let $p\in [1,\infty)$. Let $G$ be a countable discrete Abelian group which is not trivial, and $A$ be a unital separable $L^p$-operator algebra which has unique $L^p$-operator matrix norms. Then
\begin{enumerate}
\item[(i)] $\overline{M}_{G}^{p}\otimes_p A$ is $p$-isometrically incompressible if and only if $A$ is $p$-isometrically incompressible;
\item[(ii)] if $A$ is $p$-isometrically incompressible and $G$ is not trivial, then $F^p(\hat{G},F^p(G,A,\alpha),\hat{\alpha})$ is $p$-isometrically incompressible if and only if $p=2$.
\end{enumerate}
\end{prop}

\begin{proof}
(i) ``$\Longleftarrow$". Assume that $A$ is $p$-isometrically incompressible. To prove $\overline{M}_{G}^{p}\otimes_p A$ is $p$-isometrically incompressible, by \cite[Lemma 4.2]{Gardella and Thiel}, it suffices to prove $M_n^p\otimes_p A$ is $p$-isometrically incompressible for all positive integer $n$.
Let $\rho:M_n^p\otimes_p A\rightarrow \mathcal{B}\left(L^p(X,\mu)\right)$ be a contractive injective homomorphism. It suffices to prove that $\rho$ is isometric.

Since $A$ is a separable $L^p$-operator algebra, it follows that $M_n^p\otimes_p A$ is also separable.
By a similar method in \cite[Lemma 5.3]{Gardella Cuntz} and \cite[Lemma 2.7]{Lp AF}, we may assume that $(X,\mathcal{B},\mu)$ is a $\sigma$-finite measure space and $L^p(X,\mu)$ is separable.

{\it Claim 1.} For $\sum_{i,j=1}^n e_{i,j}\otimes a_{i,j}\in M_n^p\otimes_p A$, we have $$\left\|\sum_{i,j=1}^n e_{i,j}\otimes a_{i,j}\right\|=\left\|\sum_{i,j=1}^ne_{i,j}\otimes \rho(I_{M_n^p}\otimes a_{i,j})\right\|.$$

Let $\iota_1$ be the homomorphism from $M_n^p$ to $M_n^p\otimes_p A$ sending $T$ to $T\otimes I_A$, and $\iota_2$ be the homomorphism from $A$ to $M_n^p\otimes_p A$ sending $a$ to $I_{M_n^p}\otimes a$ for $T\in M_n^p$ and $a\in A$, where $I_A$ is the unit of $A$ and $I_{M_n^p}$ is the unit of $M_n^p$. Obviously, $\rho\circ \iota_2:A\rightarrow \mathcal{B}\left(L^p(X,\mu)\right)$ is a contractive injective homomorphism. Since $A$ is $p$-isometrically incompressible, it follows that $\rho\circ \iota_2$ is isometric.
Since $A$ has unique $L^p$-operator matrix norms, we have that
$$\begin{aligned}  \left\|\sum_{i,j=1}^n e_{i,j}\otimes a_{i,j}\right\|&=
\left\|\mathrm{id}_{M_n^p}\otimes (\rho\circ\iota_2)\left(\sum_{i,j=1}^n e_{i,j}\otimes a_{i,j}\right)\right\| \\ &=
\left\|\sum_{i,j=1}^ne_{i,j}\otimes \rho(I_{M_n^p}\otimes a_{i,j})\right\|.
\end{aligned}$$
This proves Claim 1.

{\it Claim 2.} $\left\|\sum_{i,j=1}^ne_{i,j}\otimes \rho(I_{M_n^p}\otimes a_{i,j})\right\|\leq \left\|\rho\left(\sum_{i,j=1}^n e_{i,j}\otimes a_{i,j}\right)\right\|.$

For convenience, we just give the proof in the case that $n=3$. The proof for the general case is similar.
Let $V_1$ be the $9\times 9$ matrix by interchanging the second row and the $4$-th row of $I_9$, $V_2$ be the $9\times 9$ matrix by interchanging the third row and the $7$-th row of $I_9$, $V_3$ be the $9\times 9$ matrix by interchanging the $6$-th row and the $8$-th row of $I_9$, where $I_9$ is the unit matrix of order $9$. One can check that $V_i^2=I_{9}$ and $V_iV_j=V_jV_i$ for all $i,j\in\{1,2,3\}$.

Since $M_3^p$ is $p$-isometrically incompressible, it follows that $\rho\circ\iota_1:M_3^p\rightarrow \mathcal{B}\left(L^p(X,\mu)\right)$ is an isometric homomorphism.
Since $M_3^p$ has unique $L^p$-operator matrix norms, it follows that $$\left\|(\mathrm{id}_{M_3^p}\otimes\left(\rho\circ \iota_1)\right)(V_i)\right\|=\left\|V_i\right\|=1$$ for all $i\in\{1,2,3\}$. Let $W_i=\left(\mathrm{id}_{M_3^p}\otimes(\rho\circ \iota_1)\right)(V_i)$. One can check that

$$W_{3} W_2W_1\left(\sum_{i,j=1}^3e_{i,j}\otimes \rho(I_{M_3^p}\otimes a_{i,j})\right)W_1W_2W_{3}=I_{M_3^p}\otimes \rho\left(\sum_{i,j=1}^3 e_{i,j}\otimes a_{i,j}\right).$$

%$$\begin{aligned} & W_{n-1}\cdots W_2W_1(\sum_{i,j=1}^ne_{i,j}\otimes \rho(I_{M_n^p}\otimes a_{i,j}))W_1W_2\cdots W_{n-1}\\
%=&
%I_{M_n^p}\otimes \rho(\sum_{i,j=1}^n e_{i,j}\otimes a_{i,j}).
%\end{aligned}$$

Hence we have
$$\begin{aligned} & \left\|\sum_{i,j=1}^3e_{i,j}\otimes \rho\left(I_{M_3^p}\otimes a_{i,j}\right)\right\|\\
=&
\left\|W_3^2W_2^2W_1^2 \left(\sum_{i,j=1}^3e_{i,j}\otimes \rho\left(I_{M_3^p}\otimes a_{i,j}\right)\right)W_1^2W_2^2W_3^2\right\| \\
 \leq &
\left\|W_3W_2W_1\left(\sum_{i,j=1}^3e_{i,j}\otimes \rho\left(I_{M_3^p}\otimes a_{i,j}\right)\right)W_1W_2W_3\right\|\\=&
\left\| I_{M_3^p}\otimes \rho\left(\sum_{i,j=1}^3 e_{i,j}\otimes a_{i,j}\right)\right\|\\=&
\left\|\rho\left(\sum_{i,j=1}^3 e_{i,j}\otimes a_{i,j}\right)\right\|.
\end{aligned}$$
This proves Claim 2.

Combing Claim 1 and Claim 2, we obtain $$\left\|\rho\left(\sum_{i,j=1}^n e_{i,j}\otimes a_{i,j}\right)\right\|=\left\|\sum_{i,j=1}^n e_{i,j}\otimes a_{i,j}\right\|.$$
Hence $\rho$ is isometric, which implies that $M_n^p\otimes_p A$ is $p$-isometrically incompressible.

``$\Longrightarrow$". Now we assume that $\overline{M}_{G}^{p}\otimes_p A$ is $p$-isometrically incompressible.
To prove $A$ is $p$-isometrically incompressible, we firstly prove the following useful fact.

{\it Claim 3.} We denoted by $\mathrm{Rep}_p(A)$ the class consisting of all contractive representations of $A$ on $\sigma$-finite separable $L^p$-spaces.
For each $T\in \overline{M}_{G}^{p}\otimes_p A$, we have $$\|T\|=\sup\{\|\mathrm{id}_{\overline{M}_{G}^{p}}\otimes \pi(T)\|:\pi\in \mathrm{Rep}_p(A)\}.$$

For each $T\in \overline{M}_{G}^{p}\otimes_p A$, we let
$$\|T\|_1=\sup\{\|\mathrm{id}_{\overline{M}_{G}^{p}}\otimes \pi(T)\|:\pi\in \mathrm{Rep}_p(A)\}.$$
Then for any positive integer $n$, there exists a separable $\sigma$-finite $L^p$-space $E_n$ and a contractive representation $\pi_n:A\rightarrow \mathcal{B}(E_n)$ such that $$\|\mathrm{id}_{\overline{M}_{G}^{p}}\otimes \pi_n(T)\|> \|T\|_1-\frac{1}{n}.$$
Let $\pi_0:A\rightarrow \mathcal{B}(E_0)$ be an isometric representation of $A$ on separable $\sigma$-finite $L^p$-space $E_0$. Then the $L^p$-direct sum (\cite[Definition 1.23]{N. C. Phillips Lp}) of $\{\pi_n\}_{n=0}^\infty$ is an isometric representation of $A$, which is denoted by $\pi$. One can check that $$\|T\|_1=\|\mathrm{id}_{\overline{M}_{G}^{p}}\otimes \pi(T)\|.$$
Since $A$ has unique $L^p$-matrix norms, it follows from \cite[Lemma 2.2]{WZ} that $\|T\|=\|\mathrm{id}_{\overline{M}_{G}^{p}}\otimes \pi(T)\|.$ This prove the Claim 3.

Let $\pi:A\rightarrow \mathcal{B}(E)$ be a contractive representation of $A$ on separable $\sigma$-finite $L^p$-space $E$. Then, by Claim 3,
$$\mathrm{id}_{\overline{M}_{G}^{p}}\otimes \pi:\overline{M}_{G}^{p}\otimes_p A\rightarrow \mathcal{B}(l^p(G)\otimes_p E)$$ is contractive.
Since $\overline{M}_{G}^{p}\otimes_p A$ is $p$-isometrically incompressible, it follows that $\mathrm{id}_{\overline{M}_{G}^{p}}\otimes \pi$ is isometric.
Let $\delta_t$ be the canonical basis of $l^p(G)$ and let $\delta_{e,e}\in \mathcal{B}(l^p(G))$ is the rank-one projection onto the span of $\{\delta_{e}\}$, where $e$ is the unit of group $G$.
Then for each $a\in A$,  the following holds $$\|\pi(a)\|=\|\mathrm{id}_{\overline{M}_{G}^{p}}\otimes \pi(\delta_{e,e}\otimes a)\|=\|\delta_{e,e}\otimes a)\|=\|a||.$$
This shows that $\pi$ is isometric, which implies that $A$ is $p$-isometrically incompressible.

(ii) When $p=2$, it follows from Theorem \ref{WZ} and Proposition \ref{P:key} that $F^{p}(\hat{G},F^{p}(G,A,\alpha),\hat{\alpha})$ is $p$-isometrically incompressible. Now we assume that $p\neq 2$. By Theorem \ref{WZ}, it follows that $F^{p}(\hat{G},F^{p}(G,A,\alpha),\hat{\alpha})$ is isometrically isomorphic to $F^{p}(G,F^{p}(\hat{G},A,\beta),\hat{\beta}\otimes \alpha)$. Hence, to show $F^p(\hat{G},F^p(G,A,\alpha),\hat{\alpha})$ is not $p$-isometrically incompressible, it suffices to show that $F^{p}(G,
F^{p}(\hat{G},A,\beta),\hat{\beta}\otimes \alpha)$ is not $p$-isometrically incompressible.

{\it Case 1.} $G$ is an infinite group.

By Theorem \ref{WZ}, it follows that $$\Phi_2:F^{p}\left(G,
F^{p}\left(\hat{G},A,\beta\right),\hat{\beta}\otimes \alpha\right)\rightarrow F^{p}\left(G,C_{0}(G,A),\mathrm{lt}\otimes\alpha\right)$$ is an injective contractive homomorphism, and $\Phi_2$ is not an isomorphism. Since $\Phi_2$ has dense range and is not surjective as noted in Remark \ref{2} , it follows that $\Phi_2$ is not an isometry.
We assume that $A\subset \mathcal{B}(E)$ is a closed subalgebra, where $E$ is a separable $L^p$-space. So
$$\Phi_4\circ\Phi_3\circ\Phi_2:F^{p}\left(G,
F^{p}\left(\hat{G},A,\beta\right),\hat{\beta}\otimes \alpha\right) \rightarrow \mathcal{B}\left(l^p(G)\otimes_p E\right)$$ is an injective contractive homomorphism which can not be isometric.
Hence $F^{p}(G, F^{p}(\hat{G},A,\beta),\hat{\beta}\otimes \alpha)$  is not $p$-isometrically incompressible.

{\it Case 2.} $G$ is a finite group but not a trivial group.

Since $G$ is finite, it follows that $\hat{G}$ is also finite. Since $G$ is amenable, the full $L^p$-operator crossed product coincides with the reduced $L^p$-operator crossed product.
By \cite[Theorem 2.19]{Choi and Gardela}, it follows that $$\mathrm{core}\left(F^{p}\left(G,F^{p}\left(\hat{G},A,\beta\right),\hat{\beta}\otimes \alpha\right)\right)=\mathrm{core}\left(F^{p}\left(\hat{G},A,\beta\right)\right)=\mathrm{core}\left(A\right)$$ and $$\mathrm{core}\left(F^{p}\left(G,C(G,A),\mathrm{lt}\otimes\alpha\right)\right)=\mathrm{core}\left(C(G,A)\right)=C\left(G,\mathrm{core}(A)\right).$$ Clearly, $\mathrm{core}(A)$ is not isometrically isomorphic to $C(G,\mathrm{core}(A))$.
 Hence $$\Phi_2: F^{p}\left(G,F^{p}\left(\hat{G},A,\beta\right),\hat{\beta}\otimes \alpha\right)\rightarrow F^{p}\left(G,C(G,A),\mathrm{lt}\otimes\alpha\right)$$ is not isometric. Since $\Phi_2$ is contractive and injective, it follows that $F^{p}(G,
F^{p}(\hat{G},A,\beta),\hat{\beta}\otimes \alpha)$ is not $p$-isometrically incompressible.\end{proof}

Before proving Theorem \ref{MainThm2}, we first prove two useful lemmas.
\begin{lem}[{\cite{incompressible}}]\label{uniformly}
For each positive integer $n$, the algebra $M_n^p$ is $p$-uniformly incompressible.
\end{lem}
\begin{proof}
The proof follows from \cite{incompressible}. For the readers' convenience, we recall the proof here.

Let $\varphi:M_n^p\rightarrow \mathcal{B}(L^p(X,\mu))$ be a bounded injective homomorphism.
Let $\la \omega,\xi\ra$ be the evaluation pairing $(l_n^p)'\times l_n^p\rightarrow \mathbb{C}$. Choose any $\xi_0\in l_n^p$ and $\omega_0\in (l_n^p)'$ with $\|\xi_0\|=\|\omega_0\|=\la \omega_0, \xi_0\ra=1$. We define $a_0\in M_n^p$ by $a_0\eta=\la \omega_0,\eta\ra\xi_0$ for all $\eta\in l_n^p$.
Then $a_0$ is a rank-one idempotent and $\|a_0\|=1$. Since $a_0^2=a_0$, it follows that $\|\varphi(a_0)\|\geq 1$.

Now for any $x\in M_n^p$ with $\|x\|=1$. Let $\epsilon>0$. We choose $\xi\in l_n^p$ such that $\|\xi\|=1$ and $\|x \xi\|>1-\epsilon$. Then we choose $\rho\in (l_n^p)'$ such that $\|\rho\|=1$ and $\la\rho,x\xi\ra=\|x\xi\|$. Let $b,c\in M_n^p$ be the rank-one operators given by $$b\eta=\la \rho,\eta\ra\xi_0$$ and $$c\eta=\la \omega_0,\eta\ra\xi$$ for all $\eta\in l_n^p$. One can check that  $\|b\|=\|c\|=1$ and $bxc=\|x\xi\|a_0$. Then we have that
$$\|x\xi\| \|\varphi(a_0)\|=\|\varphi(b)\varphi(x)\varphi(c)\|\leq \|\varphi\|^2\|\varphi(x)\|.$$ Since $\|x\xi\|>1-\epsilon$ and $\|\varphi(a_0)\|\geq 1$, it follows that $$\|\varphi(x)\|\geq \|\varphi\|^{-2}\|x\xi\|\|\varphi(a_0)\|\geq (1-\epsilon)\|\varphi\|^{-2}\|x\|.$$ Therefore, $\|\varphi(x)\|\geq \|\varphi\|^{-2}\|x\|$, which implies that $M_n^p$ is $p$-uniformly incompressible.
\end{proof}

\begin{lem}\label{direct}
Let $p\in [1,\infty)$ and $G$ be a countable discrete group. Then $\overline{M}_{G}^{p}$ is $p$-uniformly incompressible.
\end{lem}
\begin{proof}
Let $R\in [1,\infty)$ and let $\rho:\overline{M}_{G}^{p}\rightarrow \mathcal{B}(L^p(X,\mu))$ be an injective homomorphism with $\|\rho\|\leq R$. Then we let $\epsilon(R)=\min\{R,1/R^2\}$.

For any $a\in \overline{M}_G^p$ with $\|a\|=1$, there exists a finite subset $F$ of $G$ and $a_F\in M_F^p$ such that $\|a-a_F\|<\epsilon(R)/(4R)$.
We denote by $|F|$ the cardinality of $F$, then the algebra $M_F^p$ is isometrically isomorphic to $M_{|F|}^p$. By Lemma \ref{uniformly}, the algebra $M_F^p$ is $p$-uniformly incompressible with lower bound $R^{-2}$ which is independent of $F$.
Hence we have $$\|\rho(a_F)\|\geq \frac{\|a_F\|}{R^2}\geq \epsilon(R)\|a_F\|.$$

Then we have the following:
$$\begin{aligned}  \|\rho(a)\|&\geq
\| \rho(a_F)\|-\|\rho(a_F-a)\| \\ &\geq
\epsilon(R)\|a_F\|-R\cdot\frac{\epsilon(R)}{4R} \\ &\geq
\epsilon(R)(\|a\|-\|a-a_F\|)-\frac{\epsilon(R)}{4}\\ &\geq
\epsilon(R)-\frac{\epsilon(R)^2}{4R}-\frac{\epsilon(R)}{4}\\ &\geq
\epsilon(R)-\frac{\epsilon(R)}{2}\\ &=\frac{\epsilon(R)}{2}.
\end{aligned}$$
This shows that $\overline{M}_{G}^{p}$ is $p$-uniformly incompressible.
\end{proof}

Theorem \ref{MainThm2} follows from the following proposition, since isomorphisms between $L^p$-operator algebras preserve the $p$-incompressibility.

\begin{prop}\label{P:key2}
Let $p\in [1,\infty)$, $G$ be a countable discrete Abelian group, and $A=M_n^p$. Then
\begin{enumerate}
\item[(i)] $\overline{M}_{G}^{p}\otimes_p A$ is $p$-incompressible;
\item[(ii)] $F^p(\hat{G},F^p(G,A,\alpha),\hat{\alpha})$ is $p$-incompressible if and only if $p=2$ or $G$ is finite.
\end{enumerate}
\end{prop}

\begin{proof}
(i) Let $A=M_n^p$.
By Lemmas \ref{uniformly} and \ref{direct}, we have that $\overline{M}_{G}^{p}\otimes_p A$ is $p$-uniformly incompressible. Hence $\overline{M}_{G}^{p}\otimes_p A$ is $p$-incompressible.

(ii) If $p=2$ or $G$ is finite, since $F^p(\hat{G},F^p(G,A,\alpha),\hat{\alpha})$ is isomorphic to $\overline{M}_{G}^{p}\otimes_p A$, it follows that $F^p(\hat{G},F^p(G,A,\alpha),\hat{\alpha})$ is $p$-incompressible.

Now assume that $p\neq 2$ and $G$ is infinite. It follows from Theorem \ref{WZ} that $F^{p}(\hat{G},F^{p}(G,A,\alpha),\hat{\alpha})$ and $F^{p}(G,F^{p}(\hat{G},A,\beta),\hat{\beta}\otimes \alpha)$ are isometrically isomorphic. Hence, to show $F^p(\hat{G},F^p(G,A,\alpha),\hat{\alpha})$ is not $p$-incompressible, it suffices to show that $F^{p}(G,
F^{p}(\hat{G},A,\beta),\hat{\beta}\otimes \alpha)$ is not $p$-incompressible.

By Theorem \ref{WZ} and Remark \ref{2}, the map $$\Phi_2:F^{p}\left(G,
F^{p}\left(\hat{G},A,\beta\right),\hat{\beta}\otimes \alpha\right)\rightarrow F^{p}\left(G,C_{0}(G,A),\mathrm{lt}\otimes\alpha\right)$$ is injective, contractive, has dense range, and is not surjective. This implies that  $\Phi_2$ is not bounded below.
Recall that $A=M_n^p=\mathcal{B}(l_n^p)$. So $$\Phi_4\circ\Phi_3\circ\Phi_2:F^{p}\left(G,
F^{p}(\hat{G},A,\beta),\hat{\beta}\otimes \alpha\right) \rightarrow \mathcal{B}\left(l^p(G)\otimes_p l_n^p\right)$$ is an injective contractive homomorphism which can not be bounded below.
Therefore $F^{p}(G, F^{p}(\hat{G},A,\beta),\hat{\beta}\otimes \alpha)$  is not $p$-incompressible.\end{proof}

%%%%%%%%%%%%%%%%%%%%%%%%%%%%%%%%%%%%%%%%%%%%%%%%%%%%%

%\section*{Acknowledgements}

%The first author is Supported by the National Natural Science Foundation of China (grant numbers 12201240), the China Postdoctoral Science Foundation (grant number 2022M711310).
%The second author is supported by the National Natural Science Foundation of China (grant number 12171195).

\subsection*{Data Availability} No datasets were generated or analysed during the current study.

\subsection*{Declarations}
The authors have no competing interests to declare that
are relevant to the content of this article.

%We are much indebted to Yongjiang Duan, Kunyu Guo, Youqing Ji, Yanyue Shi, Kai Wang and Dechao Zheng for their interest and constructive comments on this work.
%The first author was supported by National Natural Science Foundation of China (Grant
%No. 11671167).

%%%%%%%%%%%%%%%%%%%%%%%%%%%%%%%%%%%%%%%%%55


\begin{thebibliography}{9}
%\bibitem{Blackadar}
% B. Blackadar, $K$-Theory for Operator Algebras, Second edition, Mathematical Sciences Research Institute Publications, 5, Cambridge University Press, Cambridge, 1998.

\bibitem{Blecher and Phillips}
D. Blecher and N. C. Phillips, $L^{p}$ operator algebras with approximate identities I, Pacific
J. Math. 303 (2019), no. 2, pp. 401--457.

\bibitem{Bonsall}
F. F. Bonsall, A minimal property of the norm in some Banach algebras, J. London Math. Soc. 29 (1954), 156--164.

\bibitem{Braga}
B. M. Braga  and A. Vignati, On the uniform Roe algebra as a Banach algebra and embeddings of $\ell_p$ uniform Roe algebras, Bull. Lond. Math. Soc. 52 (2020), no. 5, 853--870.

%\bibitem{Brown and Ozawa}
%N. P. Brown and N. Ozawa, $C^*$-algebras and finite-dimensional approximations, Graduate
%        Studies in Mathematics, 88, American Mathematical Society, Providence, RI, 2008.

%\bibitem{Blecher and Le Merdy}
%D. P. Blecher and C. Le Merdy, Operator Algebras and their Modules-an Operator Space
%Approach, London Mathematical Society Monographs, New Series, no. 30. Oxford Science
%Publications. The Clarendon Press, Oxford University Press, Oxford, 2004.

%\bibitem{Connes}
%A. Connes, An analogue of the Thom isomorphism for crossed products of a $C^{*}$-algebra by an action of $\mathbb{R}$, Adv. Math. 39 (1981), pp. 31-55.

\bibitem{Conway}
J. B. Conway, A Course in Functional Analysis, Graduate Texts in Mathematics, 96, Springer, New York, 1990.

\bibitem{Choi and Gardela}
Y. Choi, E. Gardella, H. Thiel, Rigidity results for $L^p$-operator algebras and
applications, Adv. Math. 452 (2024), Paper No. 109747, 47 pp.

\bibitem{Chung lp Roe 1}
Y. C. Chung and K. Li, Rigidity of $l^p$ Roe-type algebras, Bull. Lond. Math. Soc. 50 (2018), no. 6, 1056--1070.

%\bibitem{Chung lp Roe}
%Y. C. Chung and K. Li, Structure and $K$-theory of $\ell^p$ uniform Roe algebras, J. Noncommut. Geom. 15 (2021), no. 2, 581--614.

\bibitem{Cort}
G. Corti\~{n}as and M. E. Rodr\'{\i}guez, $L^p$-operator algebras associated with oriented graphs, J. Operator Theory 81 (2019), no. 1, 225--254.

\bibitem{Gardella and Thiel groupoid}
E. Gardella and M. Lupini, Representations of \'{e}tale groupoids on $L^{p}$ spaces,
Adv. Math. 318 (2017), pp. 233--278.

\bibitem{Gardella and Thiel}
E. Gardella and H. Thiel, Group algebras acting on  $L^{p}$ spaces, J. Fourier Anal. Appl. 21 (2015),
no. 6, pp. 1310--1343.

\bibitem{Gardella and Thiel quotient}
E. Gardella and H. Thiel, Quotients of Banach algebras acting on $L^{p}$ spaces,
Adv. Math. 296 (2016), pp. 85--92.

\bibitem{Gardella and Thiel convolution operators}
E. Gardella, H. Thiel, Isomorphisms of algebras of convolution
operators, Ann. Sci. \'{E}c. Norm. Sup\'{e}r. (4) 55 (2022), no. 5, 1433--1471.

\bibitem{Gardella and Thiel convolution}
E. Gardella and H. Thiel, Representations of $p$-convolution algebras on
$L^{q}$-spaces. Trans. Amer. Math. Soc. 371 (2019), no. 3, pp. 2207--2236.

\bibitem{Gardella modern}
E. Gardella, A modern look at algebras of operators on $L^{p}$ spaces,
Exo. Math. 39 (2021), no. 3, pp. 420--453.


\bibitem{Gardella Cuntz}
E. Gardella and  S. Tinghammar, A Cuntz-Krieger uniqueness theorem for $L^p$-operator graph algebras, Prepint (arXiv:2502.15591) 2025.


\bibitem{Hejazian and Pooya}
S. Hejazian and S. Pooya, Simple reduced $L^p$-operator crossed products with unique trace, J. Operator Theory 74 (2015), no. 1, 133--147.

\bibitem{Herz}
C. Herz, The theory of $p$-spaces with an application to convolution operators,
Trans. Amer. Math. Soc. 154 (1971), pp. 69--82.


\bibitem{Hew}
E. Hewitt and K. A. Ross, Abstract Harmonic Analysis. Vol. II: Structure and Analysis for Compact groups. Analysis on Locally Compact Abelian Groups, Die Grundlehren der Mathematischen Wissenschaften, Band 152
Springer-Verlag, New York-Berlin, 1970.

%\bibitem{Kaniuth}
       % E. Kaniuth, A Course in Commutative Banach Algebras (Springer, New York), 2009.


%\bibitem{Johson}
%  B. E. Johnson, An introduction to the theory of centralizers, Proc. London Math. Soc. (3) 14 (1964), 299–320.

\bibitem{Katsoulis Takai}
E. Katsoulis and C. Ramsey, Crossed products of operator algebras: applications of Takai duality, J. Funct. Anal. 275 (2018), no. 5, 1173--1207.

\bibitem{Katsoulis}
E. Katsoulis and C. Ramsey, Crossed products of operator algebras, Mem. Amer. Math. Soc. 258 (2019), no. 1240, vii+85 pp.

\bibitem{Palmer}
T. W. Palmer, Banach Algebras and the General Theory of *-Algebras. Vol. I. Algebras and Banach algebras. Encyclopedia of Mathematics and its Applications, 49. Cambridge University Press, Cambridge, 1994

\bibitem{Pedersen}
G. Pedersen, $C^*$-algebras and Their Automorphism Groups, Pure Appl. Math. (Amst.)
Academic Press, London, 2018, xviii+520 pp.

%\bibitem{Piser}
%G. Pisier, Tensor products of $C^*$-algebras and operator spaces—the Connes-Kirchberg problem, London Mathematical Society Student Texts, 96, Cambridge University Press, Cambridge, 2020.

\bibitem{Phi Cuntz}
N. C. Phillips, Analogs of Cuntz algebras on $L^p$ spaces, Preprint (arXiv:1201.4196), 2012.

\bibitem{N. C. Phillips Lp}
N. C. Phillips, Crossed products of $L^{p}$ operator algebras and the $K$-theory of Cuntz algebras
on $L^{p}$ spaces, Preprint (arXiv:1309.6406), 2013.

\bibitem{PlpsOpenQues}
N. C. Phillips, "Open problems related to operator algebras on $L^{p}$ spaces", preprint,
2014, https://pdfs.semanticscholar.org/0823/5038ec45079e7721a590 21a4492da2c2b1a3.pdf.

\bibitem{Phillips look like}
N. C. Phillips, Operator algebras on $L^{p}$ spaces which "look like" $C^{*}$-algebras,
2014,
https://pages.uoregon.edu/ncp/Talks/20140527_GPOTS/LpOpAlgs_TalkSummary.pdf.

\bibitem{Lp AF}
N. C. Phillips and M. G. Viola, Classification of $L^p$ AF algebras, Internat. J. Math. 31 (2020), no. 13, 2050088, 41 pp.

\bibitem{incompressible}
N. C. Phillips, https://pages.uoregon.edu/ncp/Research/Talks/Talk_20240709_BAOA2024_Incpr/V2_Corr1/Talk_For_posting.pdf.

\bibitem{Raeburn}
I. Raeburn, On crossed products and Takai duality, Proc. Edinburgh Math. Soc. (2) 31 (1988), pp. 321-330.

%\bibitem{Rudin}
%W. Rudin, Real and complex analysis, Third edition, McGraw-Hill Book Company, New York, 1987.

\bibitem{Takai}
H. Takai, On a duality for crossed products of $C^{*}$-algebras, J. Funct. Anal. 19 (1975), pp. 25--39.

\bibitem{WW}
Q. Wang and Z. Wang, Notes on the $\ell^p$-Toeplitz algebra on $\ell^p(\mathbb{N})$, Isr. J. Math. 245
(2021), 153--163.
%
%\bibitem{Wang}
%Z. Wang and Y. Zeng, Gelfand theory of reduced group $L^{p}$ operator algebra,  Ann. Funct. Anal. 13, 14 (2022).

\bibitem{WZ}
Z. Wang and S. Zhu, On the Takai duality for $L^{p}$ operator crossed products, Math. Z. 304 (2023), no. 4, Paper No. 54, 23 pp.

\bibitem{nuclear}
Z. Wang and S. Zhu, $p$-nuclearity of $L^p$-operator crossed products. Isr. J. Math. (2025). https://doi.org/10.1007/s11856-025-2849-4.

\bibitem{W}
Z. Wang, $p$-nuclearity of reduced group  $L^p$-operator algebras, Pacific J. Math. 342 (2026), no. 2, 427--440.

\bibitem{Williams}
D. Williams, Crossed Products of $C^{*}$-Algebras, Mathematical Surveys and Monographs, 134, American Mathematical Society, Providence, RI, 2007.
\end{thebibliography}
\end{document}